\documentclass[12pt]{amsart}
\def\Int{\operatorname{Int}}
\def\diam{\operatorname{diam}}
\def\N{\mathbb{N}}

\newtheorem{Theorem}{Theorem}[section]

\newtheorem{Lemma}[Theorem]{Lemma}

\newtheorem*{Thm}{Theorem}
\theoremstyle{definition}
\newtheorem{Definition}[Theorem]{Definition}

\theoremstyle{remark}
\newtheorem*{Remark}{Remark}
\begin{document}
\sloppy
\title{Reflection groups of geodesic spaces and Coxeter groups}
\author{Tetsuya Hosaka} 
\address{Department of Mathematics, Utsunomiya University, 
Utsunomiya, 321-8505, Japan}
\date{March 14, 2004}
\email{hosaka@cc.utsunomiya-u.ac.jp}
\thanks{Partly supported by the Grant-in-Aid for Scientific Research, 
The Ministry of Education, Culture, Sports, Science and Technology, Japan, 
(No.~15740029).}
\maketitle
\begin{abstract}
H.S.M.~Coxeter showed that 
a group $\Gamma$ is a finite reflection group of an Euclidean space 
if and only if $\Gamma$ is a finite Coxeter group. 
In this paper, we define {\it reflections} of geodesic spaces in general, 
and we prove that 
$\Gamma$ is a cocompact discrete reflection group of some geodesic space 
if and only if $\Gamma$ is a Coxeter group.
\end{abstract}

\section{Introduction and preliminaries}

The purpose of this paper is to study reflection groups of geodesic spaces.
A metric space $(X,d)$ is called a {\it geodesic space} if 
for each $x,y \in X$, 
there exists an isometry $\xi:[0,d(x,y)] \rightarrow X$ such that 
$\xi(0)=x$ and $\xi(d(x,y))=y$ (such $\xi$ is called a {\it geodesic}).
We say that an isometry $r$ of a geodesic space $X$ is a 
{\it reflection} of $X$, 
if $r^2$ is the identity of $X$, 
$X\setminus F_r$ has strictly two convex connected components 
and $\Int F_r=\emptyset$, 
where $F_r$ is the fixed-point set of $r$ 
which is called the {\it wall} of $r$.
An isometry group $\Gamma$ of a geodesic space $X$ 
is called a {\it reflection group}, 
if some set of reflections of $X$ generates $\Gamma$.
Let $\Gamma$ be a reflection group of a geodesic space $X$ and 
let $R$ be the set of all reflections of $X$ in $\Gamma$.
We note that $R$ generates $\Gamma$ by definition.
Now we suppose that 
the action of $\Gamma$ on $X$ is proper, 
that is, $\{\gamma\in\Gamma\,|\, \gamma x\in B(x,N)\}$ is finite 
for each $x\in X$ and $N>0$ (cf.\ \cite[p.131]{BH}).
Then the set $\{F_r\,|\, r\in R\}$ is locally finite.
Let $C$ be a connected component of $X\setminus \bigcup_{r\in R} F_r$, 
which is called a {\it chamber}.
In Section~2, we show that 
$\Gamma C=X\setminus \bigcup_{r\in R} F_r$.
Then $\Gamma \overline{C}=X$ 
and for each $\gamma\in\Gamma$, 
either $C\cap \gamma C=\emptyset$ or $C=\gamma C$.
We say that $\Gamma$ is a {\it cocompact discrete reflection group} of $X$, 
if $\overline{C}$ is compact and $\{\gamma\in\Gamma\,|\,C=\gamma C\}=\{1\}$.
For example, 
every Coxeter group is a cocompact discrete reflection group 
of some geodesic space.

A {\it Coxeter group} is a group $W$ having a presentation
$$\langle \,S \, | \, (st)^{m(s,t)}=1 \ \text{for}\ s,t \in S \,
\rangle,$$ 
where $S$ is a finite set and 
$m:S \times S \rightarrow \N \cup \{\infty\}$ is a function 
satisfying the following conditions:
\begin{enumerate}
\item[(1)] $m(s,t)=m(t,s)$ for each $s,t \in S$,
\item[(2)] $m(s,s)=1$ for each $s \in S$, and
\item[(3)] $m(s,t) \ge 2$ for each $s,t \in S$
such that $s\neq t$.
\end{enumerate}
The pair $(W,S)$ is called a {\it Coxeter system}.
H.S.M.~Coxeter showed that 
a group $\Gamma$ is a finite reflection group of some Euclidean space 
if and only if $\Gamma$ is a finite Coxeter group. 
Every Coxeter system $(W,S)$ induces 
the Davis-Moussong complex $\Sigma(W,S)$ which is a CAT(0) space 
(\cite{D1}, \cite{D2}, \cite{M}).
Then 
the Coxeter group $W$ is a cocompact discrete reflection group 
of the CAT(0) space $\Sigma(W,S)$.

The purpose of this paper is to prove the following theorem.

\begin{Thm}
A group $\Gamma$ is a cocompact discrete reflection group 
of some geodesic space if and only if 
$\Gamma$ is a Coxeter group.
\end{Thm}

\section{Lemmas about reflection groups}

Let $\Gamma$ be a reflection group of a geodesic space $X$ and 
let $R$ be the set of all reflections of $X$ in $\Gamma$.
We note that $R$ generates $\Gamma$ and 
$\gamma r \gamma^{-1} \in R$ 
for each $r\in R$ and $\gamma\in\Gamma$.

We suppose that the action of $\Gamma$ on $X$ is proper, 
that is, $\{\gamma\in\Gamma\,|\, \gamma x\in B(x,N)\}$ is finite 
for each $x\in X$ and $N>0$.
Then the set of walls $\{F_r\,|\, r\in R\}$ is locally finite.
Indeed if there exist $x\in X$ and $N>0$ such that 
$\{F_r\,|\, F_r\cap B(x,N)\neq \emptyset\}$ is infinite, 
then 
$\{r\in R\,|\, rx\in B(x,2N)\}$ is infinite 
which means that the action of $\Gamma$ on $X$ is not proper.

Let $C$ be a chamber defined in Section~1, 
and for each reflection $r\in R$, 
let $X_r^+$ and $X_r^-$ be two convex connected components of $X\setminus F_r$, 
where $C\subset X_r^+$ and $C\cap X_r^-=\emptyset$.
We note that $r X_r^+=X_r^-$ and $r X_r^-=X_r^+$.

Here we prove some lemmas.

\begin{Lemma}\label{lem1}
Let $r\in R$, $x\in X_r^+$ and $y\in X$. 
Then $y\in X_r^-$ if and only if $d(x,ry)<d(x,y)$.
\end{Lemma}

\begin{proof}
Suppose that $x\in X_r^+$ and $y\in X_r^-$.
Let $[x,y]$ be a geodesic from $x$ to $y$ in $X$.
Then $[x,y]\cap F_r\neq\emptyset$ by the definitions of 
$F_r$, $X_r^+$ and $X_r^-$.
Let $z\in [x,y]\cap F_r$ and 
let $[x,z]$ and $[z,y]$ be the geodesics which are restrections of $[x,y]$.
Since $z\in F_r$, $rz=z$ and 
$r[z,y]$ is a geodesic from $z$ to $ry$.
Hence $[x,z]\cup r[z,y]$ is a path from $x$ to $ry$ whose length is $d(x,y)$.
On the other hand, since 
$x$ and $ry$ are in $X_r^+$ which is convex, 
every geodesic $[x,ry]$ is contained in $X_r^+$. 
Thus the path $[x,z]\cup r[z,y]$ is not geodesic 
because $[x,z]\cup r[z,y] \not\subset X_r^+$.
Hence the length of $[x,z]\cup r[z,y]$ is greater than $d(x,ry)$, 
that is, $d(x,y)>d(x,ry)$.
Thus if $x\in X_r^+$ and $y\in X_r^-$ then $d(x,ry)<d(x,y)$.

Suppose that $x\in X_r^+$ and $y\not\in X_r^-$.
Then either $y\in X_r^+$ or $y\in F_r$.
If $y\in F_r$, then $d(x,ry)=d(x,y)$.
If $y\in X_r^+$, then $ry \in X_r^-$.
Since $x\in X_r^+$ and $ry\in X_r^-$, 
$d(x,r(ry))<d(x,ry)$ by the above argument.
Hence $d(x,y)<d(x,ry)$.
Thus if $y\not\in X_r^-$ then $d(x,y)\le d(x,ry)$.
This means that if $d(x,ry)<d(x,y)$ then $y\in X_r^-$.
\end{proof}

\begin{Lemma}\label{lem2}
$C=\bigcap_{r\in R} X_r^+$.
\end{Lemma}

\begin{proof}
Since $C\subset X_r^+$ for each $r\in R$, 
$C\subset \bigcap_{r\in R} X_r^+$.

To prove $C\supset \bigcap_{r\in R} X_r^+$, 
we show that $X\setminus C\subset X\setminus \bigcap_{r\in R} X_r^+$.
Here 
$$X\setminus \bigcap_{r\in R} X_r^+=\bigcup_{r\in R}(X\setminus X_r^+) 
=\bigcup_{r\in R}(X_r^-\cup F_r)
=(\bigcup_{r\in R} X_r^-)\cup(\bigcup_{r\in R} F_r).$$
Hence we show that 
$X\setminus C\subset(\bigcup_{r\in R} X_r^-)\cup(\bigcup_{r\in R} F_r)$.
Let $x_0\in C$ and $y\in X\setminus C$.
We suppose that $y\in X\setminus \bigcup_{r\in R} F_r$.
Let $C'$ be the connected component of $X\setminus \bigcup_{r\in R} F_r$ 
such that $y\in C'$.
Here $C\neq C'$ because $y\in X\setminus C$.
Let $[x_0,y]$ be a geodesic from $x_0$ to $y$ in $X$.
Then $[x_0,y]\not\subset X\setminus \bigcup_{r\in R} F_r$ 
and $[x_0,y]\cap (\bigcup_{r\in R} F_r)\neq\emptyset$, i.e., 
$[x_0,y]\cap F_{r_0}\neq\emptyset$ for some $r_0\in R$.
Hence $[x_0,y]\not\subset X_{r_0}^+$.
Since $x_0\in X_{r_0}^+$ and $X_{r_0}^+$ is convex, 
$y\in X_{r_0}^-$, that is, 
$y\in \bigcup_{r\in R} X_r^-$.
Thus 
$X\setminus C\subset(\bigcup_{r\in R} X_r^-)\cup(\bigcup_{r\in R} F_r)$ and 
$C\supset \bigcap_{r\in R} X_r^+$.
\end{proof}

\begin{Lemma}\label{lem3}
Let $S$ be a subset of $R$ such that $C=\bigcap_{s\in S}X_s^+$.
Then $\langle S\rangle C=X\setminus \bigcup_{r\in R} F_r$, 
where $\langle S\rangle$ is the subgroup of $\Gamma$ generated by $S$.
\end{Lemma}

\begin{proof}
Since $C\subset X\setminus \bigcup_{r\in R} F_r$, 
for each $\gamma\in \langle S\rangle$, 
$$ \gamma C\subset \gamma(X\setminus \bigcup_{r\in R} F_r) 
=X\setminus \bigcup_{r\in R} \gamma F_r.$$
Here $\gamma r \gamma^{-1}$ is a reflection and 
$\gamma F_r=F_{\gamma r \gamma^{-1}}$.
Hence 
$$ \gamma C\subset X\setminus \bigcup_{r\in R} F_{\gamma r \gamma^{-1}}
=X\setminus \bigcup_{r\in R} F_r.$$
Thus $\langle S\rangle C\subset X\setminus \bigcup_{r\in R} F_r$.

We show that $X\setminus \bigcup_{r\in R} F_r \subset \langle S\rangle C$.
Let $x_0\in C$ and let $y\in X\setminus\bigcup_{r\in R} F_r$.
If $y\in C$ then $y\in C\subset \langle S\rangle C$.
Suppose that $y\not\in C$.
Since $C=\bigcap_{s\in S}X_s^+$, 
there exists $s_1\in S$ such that $y\in X_{s_1}^-$.
Then $d(x_0,s_1 y)<d(x_0,y)$ by Lemma~\ref{lem1}.
If $s_1 y\in C$ then $y\in s_1 C \subset \langle S\rangle C$.
Suppose that $s_1 y\not\in C$.
Then there exists $s_2\in S$ such that $s_1y\in X_{s_2}^-$ and 
$d(x_0,s_2s_1 y)<d(x_0,s_1y)<d(x_0,y)$ by Lemma~\ref{lem1}.
By iterating the above argument, 
we obtain a sequence 
$s_1,\dots,s_n \in S$ such that 
\begin{align*}
d(x_0,(s_n\cdots s_1)y)&<d(x_0,(s_{n-1}\cdots s_1)y)<\cdots \\
&<d(x_0,s_2s_1 y)<d(x_0,s_1y)<d(x_0,y).
\end{align*}
Since the action of $\Gamma$ on $X$ is proper, 
this sequence is finite.
Hence $(s_n\cdots s_1)y\in C$ for some $n$.
Then $y\in (s_1\cdots s_n)C\subset\langle S\rangle C$.
Thus $X\setminus \bigcup_{r\in R} F_r \subset \langle S\rangle C$.
\end{proof}

\begin{Remark}
Since $R$ generates $\Gamma$, 
$\Gamma C=X\setminus \bigcup_{r\in R} F_r$ 
by Lemmas~\ref{lem2} and \ref{lem3}, and 
$$\Gamma \overline{C}=\overline{\Gamma C}=
\overline{X\setminus\bigcup_{r\in R}F_r}=X$$ 
because $\{\gamma C\,|\,\gamma\in\Gamma\}$ 
and $\{F_r\,|\,r\in R\}$ are locally finite 
and $\Int F_r=\emptyset$ for each $r\in R$.
\end{Remark}

\section{Cocompact discrete reflection groups}

\begin{Definition}\label{def1}
A group $\Gamma$ is called a {\it cocompact discrete reflection group} 
of a geodesic space $X$, if 
\begin{enumerate}
\item[(1)] $\Gamma$ is a reflection group of $X$, 
\item[(2)] the action of $\Gamma$ on $X$ is proper, 
\item[(3)] for a chamber $C$, $\overline{C}$ is compact, and
\item[(4)] $\{\gamma\in\Gamma\,|\,C=\gamma C\}=\{1\}$.
\end{enumerate}
\end{Definition}

Let $\Gamma$ be a cocompact discrete reflection group 
of a geodesic space $X$, let $C$ be a chamber and 
let $S$ be a {\it minimal} subset of $R$ such that $C=\bigcap_{s\in S}X_s^+$ 
(i.e.\ $C\neq \bigcap_{s\in S\setminus\{s_0\}}X_s^+$ for each $s_0\in S$).
Then by Lemma~\ref{lem3}, 
$\langle S\rangle C=X\setminus \bigcup_{r\in R} F_r=\Gamma C$.
Since $\{\gamma\in\Gamma\,|\,C=\gamma C\}=\{1\}$, 
$S$ generates $\Gamma$. 
The purpose of this section is to prove that 
the pair $(\Gamma,S)$ is a Coxeter system.

We note that the set $S$ is finite.
Indeed if $N=\diam(C)$ and $x_0\in C$, then 
$S\subset \{\gamma\in\Gamma\,|\,\gamma x_0\in B(x_0,2N)\}$ 
which is finite because the action of $\Gamma$ on $X$ is proper.

\begin{Lemma}\label{lem4}
For each $s_0\in S$, 
there exists $x_0\in C$ such that 
$d(x_0,s_0x_0)<d(x_0,(\Gamma\setminus\{1,s_0\})x_0)$.
\end{Lemma}

\begin{proof}
Let $s_0\in S$ and 
let $Y_{s_0}=\bigcap_{s\in S\setminus\{s_0\}}X_s^+$. 
Then $Y_{s_0}$ is convex open subset of $X$, 
since each $X_s^+$ is convex and open.

We first show that $Y_{s_0}\cap F_{s_0}\neq\emptyset$.
Suppose that $Y_{s_0}\cap F_{s_0}=\emptyset$. 
Then $Y_{s_0}\subset X\setminus F_{s_0}=X_{s_0}^+\cup X_{s_0}^-$.
Hence $Y_{s_0}\subset X_{s_0}^+$, since $C\subset Y_{s_0}$. 
Thus 
$$Y_{s_0}=(\bigcap_{s\in S\setminus\{s_0\}}X_s^+)\cap X_{s_0}^+
=\bigcap_{s\in S}X_s^+=C.$$
This contradicts the minimality of $S$.
Hence $Y_{s_0}\cap F_{s_0}\neq\emptyset$.

Let $y_0\in Y_{s_0}\cap F_{s_0}$. 
Since $Y_{s_0}$ is open, 
there exists $\epsilon>0$ such that 
$B(y_0,\epsilon)\subset Y_{s_0}$.
Then we can take $x_0\in B(y_0,\epsilon)\cap C$ 
because $\Int F_{s_0}=\emptyset$. 

We show that 
$d(x_0,s_0x_0)<d(x_0,(\Gamma\setminus\{1,s_0\})x_0)$.
By construction, 
$$d(x_0,s_0x_0)=2 d(x_0,F_{s_0})\le 2\epsilon$$ 
and for each $s\in S\setminus\{s_0\}$, 
$$d(x_0,sx_0)=2 d(x_0,F_s)> 2\epsilon.$$
Let $\gamma\in\Gamma\setminus\{1,s_0\}$. 
By the proof of Lemma~\ref{lem3}, 
there exists a sequence $t_1,\dots,t_n\in S$ such that 
$\gamma=t_1\cdots t_n$ and 
\begin{align*}
d(x_0,x_0)&=d(x_0,(t_n\cdots t_1)\gamma x_0) \\
&<d(x_0,(t_{n-1}\cdots t_1)\gamma x_0)<\dots <d(x_0,\gamma x_0).
\end{align*}
Here $(t_{n-1}\cdots t_1)\gamma=t_n$.
If $t_n=s_0$, then $t_n\neq\gamma$ and 
$d(x_0,s_0x_0)=d(x_0,t_nx_0)<d(x_0,\gamma x_0)$.
If $t_n\in S\setminus \{s_0\}$, then 
$$ d(x_0,s_0x_0)\le 2\epsilon <d(x_0,t_nx_0)\le d(x_0,\gamma x_0).$$
Thus $d(x_0,s_0x_0)<d(x_0,\gamma x_0)$
for each $\gamma\in\Gamma\setminus\{1,s_0\}$, that is, 
$d(x_0,s_0x_0)<d(x_0,(\Gamma\setminus\{1,s_0\})x_0)$.
\end{proof}

\begin{Lemma}\label{lem5}
Let $s\in S$ and let $r\in R$.
If $sC\subset X_r^-$ then $s=r$.
\end{Lemma}

\begin{proof}
Let $s\in S$ and $r\in R$ such that $sC\subset X_r^-$.
By Lemma~\ref{lem4}, 
there exists $x_0\in C$ such that 
$d(x_0,s x_0)<d(x_0,(\Gamma\setminus\{1,s\})x_0)$.
Then $sx_0\in sC \subset X_r^-$.
By Lemma~\ref{lem1}, 
$d(x_0,rsx_0)<d(x_0,sx_0)$. 
Since 
$d(x_0,s x_0)<d(x_0,(\Gamma\setminus\{1,s\})x_0)$, 
$rs=1$.
Hence $s=r$.
\end{proof}

For each $\gamma\in\Gamma$, 
let $\ell(\gamma)$ denote the minimum length of 
word in $S$ which represents $\gamma$.
A representation $\gamma=s_1\cdots s_n$ ($s_i \in S$) is said to be 
{\it reduced}, if $\ell(\gamma)=n$.

\begin{Lemma}\label{lem6}
Let $\gamma\in\Gamma$ and let $s\in S$.
If $\ell(\gamma)\le\ell(s\gamma)$, then $\gamma C\subset X_s^+$.
\end{Lemma}

\begin{proof}
We prove this lemma by induction on the length $\ell(\gamma)$.

If $\ell(\gamma)=0$, 
then $\gamma=1$ and $\gamma C=C\subset X_s^+$.

Let $n=\ell(\gamma)$. 
Suppose that $\ell(\gamma)\le\ell(s\gamma)$.
There exists a sequence $s_1,\dots,s_n\in S$ 
such that $\gamma=s_1\cdots s_n$.
Let $\gamma'=\gamma s_n$.
Here $\gamma'=\gamma s_n=s_1\cdots s_{n-1}$ is reduced.
Then
$$ \ell(\gamma')=n-1\le \ell(s\gamma)-1 
=\ell((s\gamma s_n)s_n)-\ell(s_n)
\le \ell(s\gamma s_n)=\ell(s\gamma').$$
Hence $\ell(\gamma')\le\ell(s\gamma')$, 
where $\ell(\gamma')=n-1$.
By the inductive hypothesis, 
$\gamma'C\subset X_s^+$.
We can take $x_0\in C$ such that 
$d(x_0,s_nx_0)<d(x_0,(\Gamma\setminus\{1,s_n\})x_0)$
by Lemma~\ref{lem4}.

Here we show that $(\gamma')^{-1}s \gamma \not\in\{1,s_n\}$.
If $(\gamma')^{-1}s \gamma=1$, then $s\gamma=\gamma'$ 
which contradicts that $\ell(s\gamma)\ge n$ and $\ell(\gamma')=n-1$.
If $(\gamma')^{-1}s \gamma=s_n$, then $s\gamma=\gamma's_n=\gamma$ 
which contradicts that $s\neq 1$.
Thus $(\gamma')^{-1}s \gamma \not\in\{1,s_n\}$.

Then 
\begin{align*}
d(\gamma'x_0,\gamma x_0)&=d(x_0,s_n x_0)
<d(x_0,(\Gamma\setminus\{1,s_n\})x_0) \\
&\le d(x_0,(\gamma')^{-1}s \gamma x_0)=d(\gamma'x_0,s\gamma x_0).
\end{align*}
Hence $d(\gamma'x_0,\gamma x_0)<d(\gamma'x_0,s\gamma x_0)$.
Here $\gamma'x_0\in\gamma'C\subset X_s^+$.
By Lemma~\ref{lem1}, 
$\gamma x_0\in X_s^+$.
This means that $\gamma C \subset X_s^+$.
\end{proof}

Using above lemmas, we prove the main theorem.

\begin{Theorem}
The pair $(\Gamma,S)$ is a Coxeter system.
\end{Theorem}

\begin{proof}
It is sufficient to show the following (\cite[p.47~(F)]{Br}):
\begin{enumerate}
\item[{\bf(F)}] For each $\gamma\in \Gamma$ and $s,t\in S$ 
such that $\ell(s\gamma)=\ell(\gamma)+1$ and 
$\ell(\gamma t)=\ell(\gamma)+1$, 
either $\ell(s\gamma t)=\ell(\gamma)+2$ or $s\gamma t=\gamma$.
\end{enumerate}

Let $\gamma\in \Gamma$ and $s,t\in S$ 
such that $\ell(s\gamma)=\ell(\gamma)+1$ and 
$\ell(\gamma t)=\ell(\gamma)+1$.
There exists a reduced representation
$\gamma=s_1\cdots s_n$, where $n=\ell(\gamma)$.
Then $s\gamma=s(s_1\cdots s_n)$ and $\gamma t=(s_1\cdots s_n)t$ 
are reduced.
If $s\gamma t=s(s_1\cdots s_n)t$ is reduced, 
then $\ell(s\gamma t)=n+2=\ell(\gamma)+2$.

We suppose that 
$s\gamma t=s(s_1\cdots s_n)t$ is not reduced. 
Then $\ell(s\gamma t)\le n+1=\ell(\gamma t)$, i.e., 
$\ell(s\gamma t)\le\ell(s(s\gamma t))$.
By Lemma~\ref{lem6}, $s\gamma tC\subset X_s^+$.
Hence $\gamma tC\subset X_s^-$.
On the other hand, 
since $\ell(s\gamma)=n+1>\ell(\gamma)$, 
$\gamma C\subset X_s^+$ by Lemma~\ref{lem6}.
Thus $\gamma C\subset X_s^+$ and $\gamma tC\subset X_s^-$, 
that is, 
$C\subset\gamma^{-1}X_s^+$ and $tC\subset\gamma^{-1}X_s^-$.
Since $C\subset\gamma^{-1}X_s^+$, 
$\gamma^{-1}X_s^+=X_{\gamma^{-1}s\gamma}^+$.
Hence $\gamma^{-1}X_s^-=X_{\gamma^{-1}s\gamma}^-$.
Thus $tC\subset\gamma^{-1}X_s^-=X_{\gamma^{-1}s\gamma}^-$.
By Lemma~\ref{lem5}, 
$t=\gamma^{-1}s\gamma$.
Therefore $s\gamma t=\gamma$.
\end{proof}

Thus we obtain that 
every cocompact discrete reflection group of a geodesic space 
is a Coxeter group.
Conversely, every Coxeter group is 
a cocompact discrete reflection group of some CAT(0) space 
(\cite{D1}, \cite{D2}, \cite{M}).

\begin{Theorem}
A group $\Gamma$ is a cocompact discrete reflection group 
of some geodesic space if and only if 
$\Gamma$ is a Coxeter group.
\end{Theorem}

%

%

\begin{thebibliography}{10}
%
\bibitem {Bo}
N.~Bourbaki, 
{\it Groupes et Algebr\`{e}s de Lie}, Chapters IV-VI, 
Masson, Paris, 1981.
%
\bibitem {BH}
M.R.~Bridson and A.~Haefliger, 
{\it Metric spaces of non-positive curvature}, 
Springer-Verlag, Berlin, 1999.
%
\bibitem {Br}
K.S.~Brown, 
{\it Buildings}, Springer-Verlag, 1980.
%
\bibitem {C1}
H.S.M.~Coxeter, 
{\it Discrete groups generated by reflections}, 
Ann.\ of Math.\ {\bf 35} (1934), 588--621.
%
\bibitem {C2}
\bysame, 
{\it The complete enumeration of finite groups of the form $R_i^2=(R_iR_j)^{k_{ij}}=1$}, 
J.\ London Math.\ Soc.\ {\bf 10} (1935), 21--25.
%
\bibitem {D1}
M.W.~Davis, {\it Groups generated by reflections and aspherical 
manifolds not covered by Euclidean space}, Ann.\ of Math.\ {\bf 117} 
(1983), 293--324.
%
\bibitem {D2}
\bysame, {\it Nonpositive curvature and reflection groups},
in Handbook of geometric topology 
(Edited by R.~J.~Daverman and R.~B.~Sher), pp.\ 373--422, 
North-Holland, Amsterdam, 2002. 
%
\bibitem {Hu}
J.E.~Humphreys, 
{\it Reflection groups and Coxeter groups}, 
Cambridge University Press, 1990.
%
\bibitem {M}
G.~Moussong, 
{\it Hyperbolic Coxeter groups}, 
Ph.D.\ thesis, The Ohio State University, 1988.
%
\end{thebibliography}
\end{document}